\documentclass{amsart}

\usepackage{bm}

\theoremstyle{plain}
\newtheorem{theorem}{Theorem}[section]

\newtheorem{lemma}[theorem]{Lemma}
\newtheorem{proposition}[theorem]{Proposition}

\theoremstyle{definition}

\theoremstyle{remark}
\newtheorem{remark}[theorem]{Remark}

\newcommand{\tends}{\rightarrow}

\begin{document}

\title[Derivative at zero and finite parts of the Barnes zeta function]{Representations for the derivative at zero and finite parts of the Barnes zeta function}

\author{J. M. B. Noronha}
\address{Lus\'{i}ada University -- North (Porto), Rua Dr. Lopo de Carvalho, 4369-006 Porto, Portugal}
\email{jnoronha@por.ulusiada.pt}


\date{November 16, 2016}
 
\maketitle

\begin{abstract}
We provide new representations
for the finite parts at the poles and the derivative at zero of the Barnes zeta function in any dimension in the general case. These representations are in the forms of series and limits.
We also give an integral representation for the finite parts at the poles. Similar results are derived for an associated function, which we term homogeneous Barnes zeta function. Our expressions immediately yield analogous representations for the logarithm of the Barnes Gamma function, including the particular case also known as multiple Gamma function. 
\end{abstract}

\section{\label{intro}Introduction}

The Barnes $\zeta$-function is a multidimensional generalization of the Hurwitz $\zeta$-function. It was introduced in its general form and studied in systematic detail in an outstanding paper by Barnes \cite{Barnes1904a} a long time ago. A possible way of defining it is the following. Let $\mathcal{H}=\{ z\in\mathbb{C}:\theta -\frac{\pi}{2}< \arg z< \theta +\frac{\pi}{2}\}$ be some half-plane through the origin in the complex plane (where $\theta$ is some fixed angle). Let $a, w_1,\ldots , w_d\in\mathcal{H}$. For $\Re (\alpha)>d$, the $d$-dimensional Barnes $\zeta$-function is given by
\begin{equation}
\zeta_{{\rm B}}(\alpha,a | \bm{w})=\sum_{\bm{n}\in\mathbb{N}_0^d}(a+\bm{n}\cdot\bm{w})^{-\alpha}\; ,
\label{Barnes_zeta}
\end{equation}
where we use the notation $\mathbb{N}_0=\{0,1,2,\ldots\}$, $\bm{n}=(n_1,\ldots,n_d)$, $\bm{w}=(w_1,\ldots,w_d)$ and $\bm{n}\cdot\bm{w}=\sum_{i=1}^dn_iw_i$. The case $d=1$ and $w_1=1$ is the Hurwitz $\zeta$-function. $\zeta_{{\rm B}}$ is usually seen as a function of the complex variable $\alpha$, with $a$, $w_1$, \ldots, $w_d$ being fixed parameters. However, nothing prevents us from seeing these as variables as well. The requirement $a, w_1,\ldots , w_d\in\mathcal{H}$ stems from the need to fix a branch cut for the $\log$ function. This branch cut can fall anywhere outside $\mathcal{H}$. For all $\bm{n}\in\mathbb{N}_0^d$, we have $a+\bm{n}\cdot\bm{w}\in\mathcal{H}$, so the summand is uniquely defined.

Barnes actually defined his $\zeta$-function in terms of a contour integral of Hankel type, which extends similar existing results for the Riemann and Hurwitz $\zeta$-functions. This automatically provides the analytic continuation of (\ref{Barnes_zeta}) to $\mathbb{C}\setminus\{1, 2, \ldots, d\}$. The points 
$\alpha=1, 2, \ldots, d$ are simple poles. This representation provides an easy way of obtaining the residues of $\zeta_{{\rm B}}(\alpha,a|\bm{w})$ and its values at the non-positive integers. These are well known \cite{Barnes1904a,Kirsten:book}.
It is not useful when it comes to finding other values. In particular, the finite parts at the poles or the derivative (with respect to $\alpha$) at $\alpha=0$, $\zeta_{{\rm B}}'(0,a|\bm{w})$, cannot be obtained in this way and expressions for these quantities are scarce (the exception being the particular $\bm{w}=(1,1,\ldots,1)$ case). However, they are of special interest for several applications. 

Most notably, functional determinants on certain background manifolds are given in terms of $\zeta_{{\rm B}}'(0,a|\bm{w})$ (see \cite{DowkerKirsten2005,Choi2014bookchapter} and references therein). These functional determinants are particularly relevant in connection to analytic torsion (see e.g. \cite{Voros1987,Sarnak1987}) and in certain quantum field theory calculations (see e.g. \cite{Dowker1994a,Bordag1996,Beneventano2001,FlachiFucci2011,Dowker2014}). 
Physical applications that involve the finite parts at the poles include Bose--Einstein condensation \cite{HolthausKalinowskiKirsten1998,Noronha2016} and functional determinants in quantum field theory \cite{Dowker1994a}.

Barnes \cite{Barnes1904a} gave contour and line integral representations for $\zeta_{{\rm B}}'(0,a|\bm{w})$. Ruijsenaars \cite{Ruijsenaars2000} provided a related line integral representation for the same quantity (this author took the parameters $w_1$, \ldots , $w_d$ to be real and $\Re (a)>0$, but in fact the representation is valid for $a, w_1, \ldots , w_d\in\mathcal{H}$ with $\theta=0$). Apart from these results, explicit expressions for $\zeta_{{\rm B}}'(0,a|\bm{w})$ in the general case have not been available in the literature so far.
The `isotropic' case, $\bm{w}=(1,1,\ldots,1)$, can be treated by casting $\zeta_{{\rm B}}$ as a finite combination  of Hurwitz $\zeta$-functions (see \cite{Barnes1904b,Bordag1996,Choi1996,Adamchik2005,Adamchik2014}) or other means (see \cite{ChoiChoSrivastava2004,Young2013,Adamchik2014}). Matsumoto \cite{Matsumoto1998,Matsumoto2002} and Elizalde \cite{Elizalde2007} gave, in the 2-dimensional case, asymptotic expansions for $\zeta_{{\rm B}}'(0,a|(w_1,w_2))$ in powers of $w_2/w_1$. 
 Elizalde \cite{Elizalde2007} also gave another expression for this quantity and  generalized his results to the $d$-dimensional case, finding a recursive expression for $\zeta_{{\rm B}}'(0,a|\bm{w})$ in terms of the derivative at $\alpha=0$ and $\alpha=-1$ and the finite parts of the $(d-1)$-dimensional case. Other recursive relation, of an asymptotic nature, was obtained by Matsumoto \cite{Matsumoto2005}.
Spreafico \cite{Spreafico2009} derived  integral and series representations in the 2-dimensional case, while Beneventano and Santangelo \cite{Beneventano2001} derived, for the same case, a series representation involving incomplete $\Gamma$-functions. The procedures in these works are not generalizable, at least in any simple way, to higher dimensions. Dowker \cite{Dowker1994a} cast $\zeta_{{\rm B}}'(0,a|\bm{w})$ in terms of derivatives of the Hurwitz $\zeta$-function in the 2-dimensional case where the $w_i$ are natural numbers, using residue classes. This latter procedure allows the treatment of higher dimensional cases (Dowker \cite{Dowker1994a} did some work on the 3-dimensional case) and is, in fact, extensible to the case where the $w_i$ are rational numbers. However, the expressions soon become extremely complicated and not very helpful (see \cite{BayadBeck2014} in this connection).

As for the finite parts at the poles, apart from the contour and line integral representations of the intimately related multidimensional Barnes $\Psi$-functions (see Eq.~(\ref{relation_FP_Psi}) below), given by Barnes \cite{Barnes1904a}, the only existing results in the literature seem to be
an integral representation given in a particular 2-dimensional case by Holthaus et al. \cite{HolthausKalinowskiKirsten1998}
and, still in the 2-dimensional case, integral and series representations given by Spreafico \cite{Spreafico2009}. Again, the procedure in \cite{Spreafico2009} is not generalizable in any direct way to higher dimensions. The method of \cite{HolthausKalinowskiKirsten1998} can be applied to higher dimensions (we will use it in section~\ref{integrals}). 

In the present work, we provide new representations for the derivative at $\alpha=0$ and the finite parts at the poles of $\zeta_{{\rm B}}$ in any dimension, in the general $\bm{w}$ case. They are quite different from the existing ones. Contrary to the previous explicit results of Barnes and Ruijsenaars for general dimension and $\bm{w}$, which are in integral form, these representations are in the forms of series and limits.
As a complement, we also present an integral representation for the finite parts, which is obtainable in a straightforward way from the methods and results of \cite{HolthausKalinowskiKirsten1998,Ruijsenaars2000}. 

Moreover, we treat a related function, which we term homogeneous Barnes $\zeta$-function, $\zeta_{{\rm Bh}}$, defined by
\begin{equation}
\zeta_{{\rm Bh}}(\alpha|\bm{w})=\sum_{\bm{n}\in\mathbb{N}_0^d\setminus\{\bm{0}\}}(\bm{n}\cdot\bm{w})^{-\alpha}\; ,
\quad \Re (\alpha)>d\; .
\label{zetaBh}
\end{equation}
$\zeta_{{\rm Bh}}$ is obtained from $\zeta_{{\rm B}}$ by setting $a=0$, having previously removed the $\bm{n}=\bm{0}$ term from the series (\ref{Barnes_zeta}). The 2-dimensional instance of this function, with $w_1$, $w_2$ real, 
was studied in detail in \cite{Spreafico2009} and it was also considered in \cite{Dowker1994a,HolthausKalinowskiKirsten1998}. 
Its pole structure is the same as that of $\zeta_{{\rm B}}$. In particular, its residues are given by setting $a=0$ in the residues of $\zeta_{{\rm B}}$. 
We obtain for this function results which parallel the ones obtained for $\zeta_{{\rm B}}$. 
Additionally, we obtain the analytic continuation of $\zeta_{{\rm Bh}}$ itself in line integral form. This is required for the integral representations of the derivative and the finite parts and cannot be achieved by simply transposing the procedure or results pertaining to the $\zeta_{{\rm B}}$ case.

It is worth noting that similar representations for the derivatives of $\zeta_{{\rm B}}$ and $\zeta_{{\rm Bh}}$ at values other than $\alpha=0$, like the negative integers, are obtainable using the same procedures.

The derivative at zero and the finite parts at the poles of $\zeta_{{\rm B}}$ and $\zeta_{{\rm Bh}}$ are very related to the Barnes $\Gamma$-function, denoted here by $\Gamma_{{\rm B}}$, and to other quantities defined and studied by Barnes (we will generally use the subscript $B$ for Barnes' functions). Therefore, in order to properly set the context, we now specify some of these relations (two further relations, pertaining to $\zeta_{{\rm Bh}}$, are given at the end of section \ref{serieslimits_Bh}, in Eqs.~(\ref{relation_rho_zetaBh}) and (\ref{relation_FP_gamma_dq})). 
Barnes defined the $d+1$-th gamma modular form $\rho_{{\rm B}}(\bm{w})$, where $d$ is the dimension, by
\begin{equation}
\log \rho_{{\rm B}}(\bm{w})=-\lim_{a\rightarrow 0}\left[ \zeta_{{\rm B}}'(0,a|\bm{w})+\log a\right] \; .
\label{def_rho}
\end{equation}
(As above and throughout this work, the prime in $\zeta_{{\rm B}}'$ or $\zeta_{{\rm Bh}}'$ stands for derivative with respect to $\alpha$.) 
He also defined his generalized $\Gamma$-function by\footnote{This is the definition used by Barnes and also, for example, in \cite{Dowker1994a,Matsumoto1998,DowkerKirsten2005,Spreafico2009}. Some authors, for example Ruijsenaars \cite{Ruijsenaars2000}, use a different definition in which $\Gamma_{{\rm B}}$ differs by the factor $\rho_{{\rm B}}(\bm{w})$.}
\begin{equation}
\log\frac{\Gamma_{{\rm B}}(a|\bm{w})}{\rho_{{\rm B}}(\bm{w})}=\zeta_{{\rm B}}'(0,a|\bm{w})
\label{def_Gamma}
\end{equation}
and $\Psi$-functions by
\begin{equation*}
\Psi_{{\rm B}} ^{(n)}(a|\bm{w})=\frac{d^n}{da^n}\log \Gamma_{{\rm B}}(a|\bm{w})\, ,\;\, n=1,2,\ldots
%\label{def_Psi}
\end{equation*}
Finally, Barnes defined the $q$-th gamma modular forms, $q=1,\ldots ,d$, by
\begin{equation}
\gamma_{dq}(\bm{w})=\lim_{a\rightarrow 0}\left[ \frac{(-1)^q(q-1)!}{a^q}-\Psi_{{\rm B}} ^{(q)}(a)\right] \; .
\label{def_gamma_dq}
\end{equation}
Note that $\Gamma_{{\rm B}}$ is the natural generalization of the usual $\Gamma$-function in the present theory 
and also that $\gamma_{11}(1)$ is simply Euler's constant, $\gamma$. $\Psi_{{\rm B}} ^{(1)}$ is the generalization 
of the usual digamma function $\Psi$, while (\ref{def_Gamma}) is a generalization of the Lerch 
formula $\zeta_{{\rm H}}'(0,a)=\log(\Gamma (a)/\sqrt{2\pi})$, where $\zeta_{{\rm H}}$ is the Hurwitz $\zeta$-function. Barnes then proved the following result regarding the finite parts at the poles of $\zeta_{{\rm B}}$. For $q=1,\ldots,d$,
\begin{equation}
F.P.\left[ \zeta_{{\rm B}}(\alpha,a|\bm{w})\right]_{\alpha=q}=\frac{(-1)^q}{(q-1)!}\Psi_{{\rm B}} ^{(q)}(a|\bm{w})-H_{q-1}Res\,\left[\zeta_{{\rm B}}(\alpha ,a|\bm{w})\right]_{\alpha=q}\; ,
\label{relation_FP_Psi}
\end{equation}
where $H_k=\sum_{i=1}^k1/i$ are the harmonic numbers and $H_0=0$ ($F.P.$ and $Res$ stand for finite part and residue, respectively). Our representations below, together with Eqs.~(\ref{def_Gamma}), (\ref{relation_FP_Psi}), (\ref{relation_rho_zetaBh}) and (\ref{relation_FP_gamma_dq}), immediately yield analogous representations for $\log\Gamma_{\rm B}$, $\Psi_{{\rm B}} ^{(q)}$, $\rho_{{\rm B}}$ and $\gamma_{dq}$. The particular case of $\Gamma_{\rm B}$ with $\bm{w}=(1,1,\ldots,1)$ is also commonly known as multiple $\Gamma$-function.

The article is organized as follows. In section \ref{serieslimits_B}, we recover a series representation for $\zeta_{{\rm B}}$ given by Barnes and use it to derive series representations for the finite parts at the poles and the derivative at $\alpha=0$ of $\zeta_{{\rm B}}$. In the same section, we derive the limit representations for these quantities. As a curiosity, we apply our results to the usual $\log\Gamma$-function. In section \ref{serieslimits_Bh}, we obtain similar results for $\zeta_{{\rm Bh}}$. In order to make our results more concrete, we apply them in section \ref{serieslimits_2D} to the special case $d=2$. In section \ref{integrals_B}, from an integral representation of Ruijsenaars 
for the analytic continuation of $\zeta_{{\rm B}}$ 
\cite{Ruijsenaars2000} (given earlier in a particular case in \cite{HolthausKalinowskiKirsten1998}), we immediately give integral representations for the finite parts at the poles of $\zeta_{{\rm B}}$. Finally, in section \ref{integrals_Bh}, 
we derive the analytic continuation of $\zeta_{{\rm Bh}}$ in line integral form and, from there, 
the respective results for the finite parts and derivative at $\alpha=0$ of $\zeta_{{\rm Bh}}$ follow.

\section{Series and limit representations}

\subsection{Barnes ${\bm\zeta}$-function}

\label{serieslimits_B}

Our representations in this section are based on a series representation given by Barnes for $\zeta_{{\rm B}}$, which does not seem to have been used ever since Barnes' original work.
We now introduce it. For simplicity of writing and following Barnes' notation, we define the symbol $F[\, ]_{x=\bm{w}}$ by
\begin{multline*}
F[f(a+x)]_{x=\bm{w}}=(-1)^{d-1}\sum_{i=1}^df(a+w_i)+(-1)^{d}\sum_{i<j}f(a+w_i+w_j)\\
\mbox{}+(-1)^{d-1}\sum_{i<j<k}f(a+w_i+w_j+w_k)+\cdots +f(a+w_1+\cdots +w_d)\; ,
\end{multline*}
where $f$ is any function. The second sum is over all pairs $(i,j)$ satisfying $1\leq i<j\leq d$, the other sums being interpreted in a similar manner.

Additionally, we introduce the symbol $G[\, ]_{x=\bm{w}}$, which we define by
\begin{equation}
G[f(a+x)]_{x=\bm{w}}=(-1)^df(a)+F[f(a+x)]_{x=\bm{w}} \; .
\label{def G}
\end{equation}
We also need the Bernoullian functions $_dS_n(a)$ used by Barnes \cite{Barnes1904a}. These are extremely related to the higher order Bernoulli polynomials\footnote{N{\"{o}}rlund used the symbol $B_n^{(d)}$. We omit the superscript as we will use it shortly with a different meaning, that of derivative of order $d$.} 
$B_n(a|\bm{w})$ studied by N\"{o}rlund \cite{Norlund1922,Norlund1924} and used in much of the later literature. 
These polynomials can be defined via a generating function by\footnote{\label{Bernoulli_def}This is the definition used for example in \cite{Norlund1922,Norlund1924,Batemanproject,Kirsten:book,Elizalde2007}. Some authors, namely in \cite{Ruijsenaars2000,BayadBeck2014}, use a different convention, where these polynomials differ by a factor of $\prod_{i=1}^dw_i$.}
\begin{equation}
z^d{\rm e}^{az}\prod_{i=1}^d\frac{w_i}{{\rm e}^{w_iz}-1}=\sum_{n=0}^{\infty}B_n(a|\bm{w})\frac{z^n}{n!}\; ,\quad
|z|<\frac{2\pi}{\max \{|w_i|\}}\; .
\label{Bernoulli}
\end{equation}
$B_n(a|\bm{w})$ is a polynomial of order $n$ in $a$. Comparing with Barnes' definition of the Bernoullian functions, we have
\begin{subequations}
\label{Bernoullian}
\begin{align}
_dS_n'(a)&=\frac{n!}{(n+d-1)!\prod_{i=1}^dw_i}B_{n+d-1}(a|\bm{w})\; , \\
_dS_n(0)&=0\;,\qquad n=0,1,2,\ldots ,
\end{align}
\end{subequations}
where the prime in the first equation means derivative with respect to $a$ (the dependency on $\bm{w}$ is omitted for economy of writing).

We have the following representation for $\zeta_{{\rm B}}$ due to Barnes \cite{Barnes1904a}. 
\begin{theorem}[Barnes]
Let $k$ be an integer greater than $-d$. For $\Re(\alpha)>-k$,
\begin{multline}
\zeta_{{\rm B}}(\alpha,a|\bm{w})=\sum_{\bm{n}\in\mathbb{N}_0^d}\left[
\frac{1}{(a+\bm{n}\cdot\bm{w})^{\alpha}}\right. \\
\mbox{}\left. 
-G\left[\sum_{m=0}^{k+d-1}\frac{_dS_m^{(d)}(0)}{m!}\frac{d^m}{dx^m}\frac{(a+\bm{n}\cdot\bm{w}+x)^{d-\alpha}}{(d-\alpha)\cdots(1-\alpha)}\right]_{x=\bm{w}}\right] \\
\mbox{}+(-1)^d\sum_{m=0}^{k+d-1}\frac{_dS_m^{(d)}(0)}{m!}\frac{d^m}{da^m}\frac{a^{d-\alpha}}{(d-\alpha)\cdots(1-\alpha)} \; .
\label{Barnes_series1}
\end{multline}
For $\Re(\alpha)>-k$, the infinite series in $\bm{n}$ in (\ref{Barnes_series1}) is absolutely convergent.
\end{theorem}
The superscript $(d)$ in $_dS_m^{(d)}$ means derivative of order $d$. This representation can be cast as
\begin{multline}
\zeta_{{\rm B}}(\alpha,a|\bm{w})=\sum_{\bm{n}\in\mathbb{N}_0^d}\left[
\frac{1}{(a+\bm{n}\cdot\bm{w})^{\alpha}}
\right. \\
\left. -\sum_{m=0}^{k+d-1}\frac{_dS_m^{(d)}(0)}{m!}\frac{\Gamma(1-\alpha)}{\Gamma(d-\alpha-m+1)}
%\right.\\
%\mbox{}\left.\times
G\left[(a+\bm{n}\cdot\bm{w}+x)^{d-\alpha-m}\right]_{x=\bm{w}}
%\right. \\
%\mbox{}\left. 
\right] \\
\mbox{}
+(-1)^d\sum_{m=0}^{k+d-1}\frac{_dS_m^{(d)}(0)}{m!}\frac{\Gamma(1-\alpha)}{\Gamma(d-\alpha-m+1)}a^{d-\alpha-m}
 \; , \qquad \Re(\alpha)>-k\; .
\label{Barnes_series2}
\end{multline}

For our purposes, we need to establish some further properties of the infinite series in (\ref{Barnes_series1}) or (\ref{Barnes_series2}). First, we observe that its summand
is an entire function of $\alpha$ for all $\bm{n}\in\mathbb{N}_0^d$ and fixed $a , w_1, \ldots , w_d\in\mathcal{H}$ (the expression being interpreted in the sense of analytic continuation when $\alpha=1,\ldots,d$).
The summand is obviously analytic for $\alpha\neq 1,\ldots,d$. At $\alpha=1,\ldots,d-m$ (with $m=0,\ldots,d-1$), the factor 
$\Gamma(1-\alpha)/\Gamma(d-\alpha-m+1)$ in (\ref{Barnes_series2}) has simple poles. 
However, at these values of $\alpha$, the factor $G\left[(a+\bm{n}\cdot\bm{w}+x)^{d-\alpha-m}\right]_{x=\bm{w}}$ 
vanishes, canceling the behaviour of the $\Gamma$-function. In fact, from combinatorial considerations we have that for any constant $c$, $G[c]_{x=\bm{w}}=0$. In addition, for $n=1,\ldots,d-1$, $G[x^n]_{x=\bm{w}}=0$. Then, by applying a binomial expansion to $(c+x)^n$, we have that for $n=0,\ldots,d-1$, $G[(c+x)^n]_{x=\bm{w}}=0$, from where the assertion follows.

The proof by Barnes of the absolute convergence of the $\bm{n}$ series in (\ref{Barnes_series1}) (\cite{Barnes1904a}, pp. 393--395) can also be used to show the uniform convergence of this series in the variable $\alpha$. Barnes' proof is rather long and involved and we will not retrace all the steps here. We merely invoke the intermediate results within it that are necessary to show uniform convergence.

Barnes showed that except for a finite number of terms (the exact number of which depends only on $a$ and $\bm{w}$), the summand of the series can be cast as $\sum_{q=0}^{\infty}c_q$, where
\begin{align*}
c_q&=\frac{\Gamma(1-\alpha)}{\Gamma(1-\alpha-q-k-d)}
\frac{b_q}{(a+\bm{n}\cdot\bm{w})^{\alpha +q+k+d}}\\
\intertext{and}
b_q&=\sum_{m=0}^{k+d-1}\frac{_dS_m^{(d)}(0)}{m!(q+k+2d-m)!}F\left[ x^{q+k+2d-m}\right]_{x=\bm{w}}\; .
\end{align*}
Now,
\begin{equation*}
\frac{c_{q+1}}{c_q}=-\frac{\alpha+q+k+d}{a+\bm{n}\cdot\bm{w}}\frac{b_{q+1}}{b_q}\; .
\end{equation*}
Barnes observed that $|(\alpha+q+k+d)b_{q+1}/b_q|$ is bounded. More precisely, it is possible to find $\lambda>0$ such that, for all $q\in\mathbb{N}_0$, we have $|(\alpha+q+k+d)b_{q+1}/b_q|<\lambda$. ($\lambda$ is dependent on $\alpha$, $k$, $d$ and $\bm{w}$. However, note that $k$, $d$ and $\bm{w}$, as well as $a$, are all being taken as fixed here.) To show uniform convergence, let $A>0$ be some constant. From the above, it is clear that there is $\lambda'>0$ such that, for all $q\in\mathbb{N}_0$ and all $\alpha\in\mathbb{C}$ such that $|\alpha|<A$, we have $|(\alpha+q+k+d)b_{q+1}/b_q|<\lambda'$.

Let $\eta$ be any constant satisfying $\eta\in(0,1)$. Then, provided $|a+\bm{n}\cdot\bm{w}|\geq \lambda'/\eta$, we have $|c_{q+1}/c_q|<\eta$. Hence, except for a finite number of terms of the $\bm{n}$ series (the exact number of which being bounded provided $|\alpha|<A$), its summand, $\sum_{q=0}^{\infty}c_q$, satisfies
\begin{equation*}
\left|\sum_{q=0}^{\infty}c_q\right|\leq |c_0|\sum_{q=0}^{\infty}\eta^{q}
=
\frac{|b_0|}{1-\eta}\left|\frac{\Gamma(1-\alpha)}{\Gamma(1-\alpha-k-d)}\right|
\frac{1}{\left|(a+\bm{n}\cdot\bm{w})^{\alpha +k+d}\right|}\; .
\end{equation*}
Now, for $|\alpha|<A$, $|\Gamma(1-\alpha)/\Gamma(1-\alpha-k-d)|$ is bounded. In addition, it is clear that the series $\sum_{\bm{n}\in\mathbb{N}_0^d}\left|(a+\bm{n}\cdot\bm{w})^{-\alpha-k-d}\right|$ is uniformly convergent in $\alpha$ provided $|\alpha|<A$ and $\Re(\alpha)>-k+\delta$ (where $\delta$ is any positive constant). Therefore, we have
\begin{proposition}
The infinite series in (\ref{Barnes_series1}) is uniformly convergent in $\alpha$ in any compact domain contained in $\{\alpha\in\mathbb{C}:\Re(\alpha)>-k\}$.
\end{proposition}

It follows that the infinite series term in (\ref{Barnes_series1}) or (\ref{Barnes_series2}) is an analytic function of $\alpha$, for $\Re(\alpha)>-k$. The pole structure of $\zeta_{{\rm B}}(\alpha,a|\bm{w})$ is therefore contained in the last term of (\ref{Barnes_series1}). From this term, one can easily recover the residues of $\zeta_{{\rm B}}$. However, our aim here is to obtain the finite parts and the derivative at $\alpha=0$.

\subsubsection{Finite parts}

In order to obtain the finite part at $\alpha=q$ ($q=1,\ldots,d$), take $k=1-q$ in (\ref{Barnes_series2}). Then, the series term is an analytic function of $\alpha$ for $\Re(\alpha)>q-1$. The finite part is given by the value of the series at $\alpha=q$ plus the finite part of the term outside the series. For the value of the series part, since $G\left[(a+\bm{n}\cdot\bm{w}+x)^{d-\alpha-m}\right]_{x=\bm{w}}$ vanishes at $\alpha=q$, we must expand it around this point. The finite part of the term outside the series is easily obtained. We finally have
\begin{multline}
F.P.[\zeta_{{\rm B}}(\alpha,a|\bm{w})]_{\alpha=q}=
\frac{(-1)^{d-q+1}}{(q-1)!}\sum_{m=0}^{d-q} \frac{_dS_m^{(d)}(0)a^{d-q-m}}{m!(d-q-m)!}(\log a-H_{d-q-m}+H_{q-1})\\
\mbox{}+\sum_{\bm{n}\in\mathbb{N}_0^d}\left[
\frac{1}{(a+\bm{n}\cdot\bm{w})^{q}} +\frac{(-1)^q}{(q-1)!}\right. \\
\mbox{}\left. \times
\sum_{m=0}^{d-q}\frac{_dS_m^{(d)}(0)}{m!(d-q-m)!}G\left[ (a+\bm{n}\cdot\bm{w}+x)^{d-q-m}\log(a+\bm{n}\cdot\bm{w}+x)\right]_{x=\bm{w}}\right] \; .
\label{fp_B_series}
\end{multline}

The finite parts can also be given the form of limits. In order to do this, it is convenient to introduce some auxiliary notation and a lemma. Let $u(\bm{n})$ be any complex-valued function of $\bm{n}\in\mathbb{N}_0^d$. For any $\bm{v}\in\mathbb{N}_0^d$, let $\bm{v}_{(i_1,i_2,\ldots,i_k)}$ (with $1\leq i_1<i_2<\cdots <i_k\leq d$) be the vector obtained from $\bm{v}$ by making all coordinates zero except the $i_1,i_2,\ldots,i_k$ coordinates. So, for example, $\bm{v}_{(2)}=(0,v_2,0,\ldots,0)$ and $\bm{v}_{(1,d)}=(v_1,0,\ldots,0,v_d)$. Define the symbol $[u(\bm{n})]_{\bm{v}}$ by
\begin{multline*}
[u(\bm{n})]_{\bm{v}}=(-1)^du(\bm{n})+(-1)^{d+1}\sum_{i=1}^du\left( \bm{n}+\bm{v}_{(i)}\right)\\
\mbox{}
+(-1)^d\sum_{1\leq i<j\leq d}u\left( \bm{n}+\bm{v}_{(i,j)}\right)+\cdots+u(\bm{n}+\bm{v})\; .
\end{multline*}
Define also $\mathcal{C}_M=\{0,1,\ldots,M\}^d$. We wish to consider sums of the type $\sum_{\bm{n}\in\mathcal{C}_M}[u(\bm{n})]_{\bm{1}}$, where $\bm{1}=(1,1,\ldots,1)$ and $M$ is any positive integer. The following lemma shows that these are multidimensional analogues of telescoping sums. As such, all of its terms cancel except the ones at the edges, yielding a simple result.
\begin{lemma}
Being $[u(\bm{n})]_{\bm{v}}$ defined as above and $M$ any non-negative integer, we have
\begin{equation*}
\sum_{\bm{n}\in\mathcal{C}_M}[u(\bm{n})]_{\bm{1}}=[u(\bm{0})]_{(M+1)\bm{1}}\; .
\end{equation*}
\end{lemma}
\begin{remark}
This is a generalization to any dimension of the trivial property of telescoping sums $\sum_{n=0}^M(u(n+1)-u(n))=u(M+1)-u(0)$.
\end{remark}
\begin{proof}[Proof of the Lemma]
We start by grouping the terms of the sum in hyper-cubic shells. Let the $k$-th shell be $S_k=\{\bm{n}\in\mathbb{N}_0^d: \max\{n_1,\ldots,n_d\} =k\}$. We now find the sum of the terms in the $k$-th shell:
\begin{equation}
\sum_{\bm{n}\in S_k}[u(\bm{n})]_{\bm{1}}\; .
\label{k-shell sum}
\end{equation}
For $k\geq 1$, the elements of $S_k\cup S_{k+1}$ are of one of the four following types:
\begin{subequations}
\label{types}
\begin{align}
i)\; & \bm{p}=(\underbrace{k+1,\ldots,k+1}_r,\underbrace{0,\ldots,0}_s,\underbrace{p_{r+s+1},\ldots,p_d}_{d-r-s})\, ,&
1\leq p_i\leq k \text{ for }i>r+s\label{type1}\\
ii)\; &\bm{p}=(\underbrace{k+1,\ldots,k+1}_r,\underbrace{0,\ldots,0}_s)\, ,& r+s=d\label{type2}\\
iii)\;  & \bm{p}=(\underbrace{k,\ldots,k}_r,\underbrace{0,\ldots,0}_s,\underbrace{p_{r+s+1},\ldots,p_d}_{d-r-s})\, ,
&1\leq p_i\leq k-1 \text{ for }i>r+s\label{type3}\\
iv)\; &\bm{p}=(\underbrace{k,\ldots,k}_r,\underbrace{0,\ldots,0}_s)\, ,& r+s=d\label{type4}
\end{align}
\end{subequations}
or elements obtained from these ones by permutation of the $p_i$'s, $i=1,\ldots,d$. In all four types, we must have $r\geq 1$ and, in the first and third types, $d-r-s\geq 1$ ($s$ is allowed to be zero). The third type does not occur if $k=1$.

Only terms of the form $u(\bm{p})$, with $\bm{p}$ being as above, contribute to the sum (\ref{k-shell sum}). We will consider one at a time the contribution of the respective four types of terms to this sum. 
We will take $\bm{p}$ as being in one of the forms displayed in (\ref{types}), $1\leq r\leq d$. Naturally, the contributions from terms $u(\bm{p})$, where $\bm{p}$ is obtained from permutation of the coordinates in (\ref{types}), will be analogous.
From the definition of $[u(\bm{n})]_{\bm{1}}$, we see that a particular term $u(\bm{p})$, with $\bm{p}$ as in (\ref{type1}), appears in the sum (\ref{k-shell sum}):
\begin{itemize}
\item once, with sign $(-1)^{d+r}$, coming from $\bm{n}=(k,\ldots,k,0,\ldots,0,p_{r+s+1},\ldots,p_d)$;
\item $\binom{d-r-s}{1}$ times, with sign $(-1)^{d+r+1}$, coming from\\
$\bm{n}=(k,\ldots,k,0,\ldots,0,p_{r+s+1},\ldots,p_d)-\bm{1}_{(i)}$, $i=r+s+1,\ldots,d$;
\item $\binom{d-r-s}{2}$ times, with sign $(-1)^{d+r+2}$, coming from\\
$\bm{n}=(k,\ldots,k,0,\ldots,0,p_{r+s+1},\ldots,p_d)-\bm{1}_{(i,j)}$, $r+s< i<j\leq d$;
\end{itemize}
and so on. Thus, the total contribution of this term to the sum is\\ $\sum_{i=0}^{d-r-s}(-1)^{d+r+i}\binom{d-r-s}{i}u(\bm{p})=0$. Hence, terms of the first type totally cancel out in (\ref{k-shell sum}). 

As for terms of the second type, the sole contribution to the appearance of the term $u(k+1,\ldots,k+1,0,\ldots,0)$ in (\ref{k-shell sum}) is from $\bm{n}=(k,\ldots,k,0,\ldots,0)$, with sign $(-1)^{d+r}$. The contribution is thus $(-1)^{d+r}u(\bm{p})$, with $\bm{p}$ as in (\ref{type2}).

Consider now a particular term of the third type (applicable only if $k>1$): $u(\bm{p})$ with $\bm{p}$ as in (\ref{type3}). On one hand, contributions to this term arise when $\bm{n}$ in (\ref{k-shell sum}) is of one of the forms $\bm{n}=\bm{p}$, $\bm{n}=\bm{p}-\bm{1}_{(i)}$ with $r+s<i\leq d$, $\bm{n}=\bm{p}-\bm{1}_{(i,j)}$ with $r+s<i<j\leq d$, \ldots, $\bm{n}=\bm{p}-\bm{1}_{(r+s+1,\ldots,d)}$. The total from these contributions is obtained as before as $\sum_{i=0}^{d-r-s}(-1)^{d+i}\binom{d-r-s}{i}u(\bm{p})=0$. Lowering one or more (but not all) of the $n_i$'s, $i=1,\ldots,r$,  in these $\bm{n}$ from $k$ to $k-1$, we obtain analogous sets of contributions to the appearance of $u(\bm{p})$ in (\ref{k-shell sum}). Likewise, all these contributions vanish. Thus, all terms of the third type cancel out in (\ref{k-shell sum}).

Finally, we consider terms of the fourth type: $u(\bm{p})$ with at least one $p_i=k$, the other $p_i$'s (possibly none) being zero. Let $\bm{p}$ be as in (\ref{type4}). There is a contribution to the appearance of such term in (\ref{k-shell sum}) coming from $\bm{n}=\bm{p}$; it is $(-1)^du(\bm{p})$. A second 
contribution comes from lowering one of the $k$'s to $k-1$, i.e., $\bm{n}=\bm{p}-\bm{1}_{(i)}$, 
$i=1,\ldots,r$. 
This contribution is $(-1)^{d+1}\binom{r}{1}u(\bm{p})$. Carrying on, we obtain contributions from lowering 
by 1 any number of $k$'s up to $r-1$. The total contribution is thus $\sum_{i=0}^{r-1}(-1)^{d+i}\binom{r}{i}u(\bm{p})=(-1)^{d+r+1}u(\bm{p})$.

Summing up, for $k\geq 1$ we have
\begin{multline}
\sum_{\bm{n}\in S_k}\left[u(\bm{n})\right]_{\bm{1}}=
\left[u\left((k+1)\bm{1}\right)-u\left(k\bm{1}\right)\right]-
\sum_{1\leq i_1<\cdots <i_{d-1}\leq d}\left[u\left((k+1)\bm{1}_{(i_1,\ldots,i_{d-1})}\right)\right.\\
\mbox{}\left. -u\left(k\bm{1}_{(i_1,\ldots,i_{d-1})}\right)\right]
+\cdots
+(-1)^{d}\sum_{1\leq i_1<i_{2}\leq d}\left[ u\left((k+1)\bm{1}_{(i_1,i_2)}\right)-u\left(k\bm{1}_{(i_1,i_2)}\right)\right] \\
\mbox{}+(-1)^{d+1}\sum_{1\leq i\leq d}\left[u\left((k+1)\bm{1}_{(i)}\right)-u\left(k\bm{1}_{(i)}\right)\right] \; .
\label{k-shell sum extended}
\end{multline}
The terms where $k+1$ appear are the contributions of the second type, while the terms where $k$ appears are the contributions of the fourth type. The expression in the first square brackets comes from the $r=d$ case, the one in the second square brackets comes from the $r=d-1$ case, \textit{etc}. By direct inspection, we see that this equality is satisfied also in the $k=0$ case. Moreover, using the notation introduced before the lemma, Eq.~(\ref{k-shell sum extended}) can be written in compact form. Thus, for all $k$,
\begin{equation}
\sum_{\bm{n}\in S_k}\left[u(\bm{n})\right]_{\bm{1}}=\left[u(\bm{0})\right]_{(k+1)\bm{1}}-\left[u(\bm{0})\right]_{k\bm{1}}\; .
\label{k-shell sum compact}
\end{equation}
(Note that $[u(\bm{0})]_{\bm{0}}=0$.) 

Using (\ref{k-shell sum compact}) in $\sum_{\bm{n}\in\mathcal{C}_M}[u(\bm{n})]_{\bm{1}}=
\sum_{k=0}^M\sum_{\bm{n}\in S_k}[u(\bm{n})]_{\bm{1}}$ and noting that the result is a telescoping sum, the lemma follows at once.
\end{proof}
We now return to the finite parts of $\zeta_{{\rm B}}$, given in (\ref{fp_B_series}). These can be cast as
\begin{multline}
F.P.[\zeta_{{\rm B}}(\alpha,a|\bm{w})]_{\alpha=q}=\lim_{M\tends\infty}\left[
\sum_{\bm{n}\in\mathcal{C}_{M-1}}(a+\bm{n}\cdot\bm{w})^{-q}
+\frac{(-1)^q}{(q-1)!}
\right.\\
\mbox{}\left.\times
\sum_{m=0}^{d-q}
\frac{_dS_m^{(d)}(0)}{m!(d-q-m)!}
\sum_{\bm{n}\in\mathcal{C}_{M-1}}
G\left[ (a+\bm{n}\cdot\bm{w}+x)^{d-q-m}\log(a+\bm{n}\cdot\bm{w}+x)\right]_{x=\bm{w}}\right] \\
\mbox{}+\frac{(-1)^{d-q+1}}{(q-1)!}\sum_{m=0}^{d-q} \frac{_dS_m^{(d)}(0)a^{d-q-m}}{m!(d-q-m)!}(\log a-H_{d-q-m}+H_{q-1})\; .
\label{fp_B_lim1}
\end{multline}
Consider the summand in the last sum inside the curly brackets. By choosing $u(\bm{n})=
(a+\bm{n}\cdot\bm{w})^{d-q-m}\log(a+\bm{n}\cdot\bm{w})$, we have
\[
G\left[ (a+\bm{n}\cdot\bm{w}+x)^{d-q-m}\log(a+\bm{n}\cdot\bm{w}+x)\right]_{x=\bm{w}}=\left[ u(\bm{n})\right]_{\bm{1}}\; .
\]
Then, upon application of the lemma to Eq.~(\ref{fp_B_lim1}), noting in addition that $\left[ u(\bm{0})\right]_{M\bm{1}}=G\left[ (a+\bm{n}\cdot\bm{w}+x)^{d-q-m}\log(a+\bm{n}\cdot\bm{w}+x)\right]_{x=M\bm{w}}$ and using Eq.~(\ref{def G}), we obtain
\begin{multline}
F.P.[\zeta_{{\rm B}}(\alpha,a|\bm{w})]_{\alpha=q}=%\\
\lim_{M\tends\infty}\Biggl[
\frac{(-1)^q}{(q-1)!}%\right.
\\
\mbox{}\left.\times
\sum_{m=0}^{d-q}\frac{_dS_m^{(d)}(0)}{m!(d-q-m)!}
F\left[ (a+x)^{d-q-m}\log(a+x)\right]_{x=M\bm{w}}\right.
\\
\left. 
+\sum_{\bm{n}\in\mathcal{C}_{M-1}}(a+\bm{n}\cdot\bm{w})^{-q}
\right]
%\\
%\mbox{}
+\frac{(-1)^{d-q+1}}{(q-1)!}\sum_{m=0}^{d-q} \frac{_dS_m^{(d)}(0)a^{d-q-m}}{m!(d-q-m)!}(H_{q-1}-H_{d-q-m})\; .
\label{fp_B_lim2}
\end{multline}

\subsubsection{Derivative at $\alpha=0$}

For the derivative at $\alpha=0$, we take $k=1$ in (\ref{Barnes_series2}), so that the representation is valid for $\Re(\alpha)>-1$. At $\alpha=0$ there are no singularities involved and all we need to do is differentiate. Since the series in (\ref{Barnes_series2}) is uniformly convergent, we can perform the differentiation directly inside the series. The result can be further simplified by making use of the property
$G[(c+x)^n]_{x=\bm{w}}=0$, $n=0,1,\ldots,d-1$ and also $G[(c+x)^d]_{x=\bm{w}}=d!\prod_{i=1}^dw_i$ (which again can be seen from combinatorial considerations). In addition, we use from \cite{Barnes1904a}, $_dS_0^{(d)}(0)=1/\prod_{i=1}^dw_i$. We thus obtain
\begin{multline}
\zeta_{{\rm B}}'(0,a|\bm{w})
=(-1)^{d+1}\sum_{m=0}^{d} \frac{_dS_m^{(d)}(0)a^{d-m}}{m!(d-m)!}(\log a-H_{d-m})
+\sum_{\bm{n}\in\mathbb{N}_0^d}\Biggl[
-\log(a+\bm{n}\cdot\bm{w}) 
\\
\mbox{}\left.
-H_d +\sum_{m=0}^{d}\frac{_dS_m^{(d)}(0)}{m!(d-m)!}G\left[ (a+\bm{n}\cdot\bm{w}+x)^{d-m}\log(a+\bm{n}\cdot\bm{w}+x)\right]_{x=\bm{w}}\right] 
\; .
\label{derivative_B_series}
\end{multline}
Similarly to what was done in the case of the finite parts, $\zeta_{{\rm B}}'(0,a|\bm{w})$ can be given in the form of a limit as
\begin{multline}
\zeta_{{\rm B}}'(0,a|\bm{w})=\lim_{M\tends\infty}\left[ -H_dM^d
+\sum_{m=0}^{d}\frac{_dS_m^{(d)}(0)}{m!(d-m)!}F\left[ (a+x)^{d-m}\log(a+x)\right]_{x=M\bm{w}}
\right.\\
\left.\mbox{}-\sum_{\bm{n}\in\mathcal{C}_{M-1}}
\log(a+\bm{n}\cdot\bm{w})\right]
+(-1)^d\sum_{m=0}^{d} \frac{_dS_m^{(d)}(0)H_{d-m}}{m!(d-m)!}a^{d-m}\; .
\label{derivative_B_lim}
\end{multline}

In the $d=1$ case, the above formulae, as well as the ones derived over the next sections, mostly reduce to well-known relations involving the Riemann and Hurwitz $\zeta$-functions and the $\Gamma$- and $\Psi$-functions. Perhaps one exception worth mentioning is the $d=1$ instance of (\ref{derivative_B_series}) and (\ref{derivative_B_lim}). As already mentioned, $\zeta_{{\rm B}}$ reduces to $\zeta_{{\rm H}}$ in this case. Then, using the Lerch formula for $\zeta_{{\rm H}}'(0,a)$, (\ref{derivative_B_series}) and (\ref{derivative_B_lim}) yield the $\log\Gamma$ representations
\begin{equation}
\log\Gamma (a)=a(\log a -1)+\frac{1}{2}\log\frac{2\pi}{a}+
\sum_{n=0}^{\infty}\left[ \left( a+n+\frac{1}{2}\right) \log\left( 1+\frac{1}{a+n}\right) -1\right]
\label{Gamma_series}
\end{equation}
and
\begin{equation*}
\log\Gamma (a)=-a+\lim_{M\rightarrow\infty}\left[
-M+\left(a+M-\frac{1}{2}\right) \log (a+M) -\sum_{n=0}^{M-1}\log (a+n)\right] \; .
\end{equation*}
Another relation can be obtained, for $|a|>1$, $\Re (a)>0$, by Taylor expanding the logarithm in (\ref{Gamma_series}) in powers of $1/(a+n)$ and interchanging the summation signs of the resulting double sum (allowed since it is absolutely convergent). This procedure yields
\begin{equation}
\log\Gamma (a)=a(\log a-1)+\frac{1}{2}\log\frac{2\pi}{a}+\sum_{k=2}^{\infty}\frac{(-1)^k(k-1)}{2k(k+1)}\zeta_{{\rm H}}(k,a)\; .
\label{Gamma_series2}
\end{equation}
We can drop the requirement $|a|>1$ of (\ref{Gamma_series2}) if we split off the $n=0$ term before expanding the logarithm, in which case we obtain a variant of this relation.

\subsection{Homogeneous Barnes ${\bm\zeta}$-function}

\label{serieslimits_Bh}

The homogeneous Barnes $\zeta$-function is defined in (\ref{zetaBh}). Since both (\ref{Barnes_zeta}) and (\ref{Barnes_series2}), with the $1/a^\alpha$ term removed from the series, are analytic functions of $a$ for $a\in\mathcal{H}$ plus a narrow strip, parallel to the border of the half-plane $\mathcal{H}$, which includes $a=0$ in its interior, the analytic continuation of $\zeta_{{\rm Bh}}(\alpha|\bm{w})$ to $\Re(\alpha)>-k$ is obtained by simply setting $a=0$ in (\ref{Barnes_series2}) (with the $1/a^\alpha$ term removed). Note that the term outside the series is cancelled by an opposite contribution coming from the series part, so that we have
\begin{multline*}
\zeta_{{\rm Bh}}(\alpha|\bm{w})=
-\sum_{m=0}^{k+d-1}\frac{_dS_m^{(d)}(0)}{m!}\frac{\Gamma(1-\alpha)}{\Gamma(d-\alpha-m+1)}
F\left[x^{d-\alpha-m}\right]_{x=\bm{w}}\\
\mbox{}
+\sum_{\bm{n}\in\mathbb{N}_0^d\setminus\{\bm{0}\}}\Biggl[
%\left[
\left.
\frac{1}{(\bm{n}\cdot\bm{w})^{\alpha}}
-\sum_{m=0}^{k+d-1}\frac{_dS_m^{(d)}(0)}{m!}\frac{\Gamma(1-\alpha)}{\Gamma(d-\alpha-m+1)}G\left[(\bm{n}\cdot\bm{w}+x)^{d-\alpha-m}\right]_{x=\bm{w}}\right]
 \, , 
 %\qquad 
 \\ \Re(\alpha)>-k\; .
%\label{hBarnes_series}
\end{multline*}
The pole structure is contained in the second term. The first term (the infinite series) is an analytic function of $\alpha$ for $\Re(\alpha)>-k$.

The procedure for obtaining the finite parts at the poles and the derivative at $\alpha=0$ is the same as for $\zeta_{{\rm B}}$. The finite parts at $\alpha=q$, $q=1,\ldots,d$, are given in the series form as
\begin{multline*}
F.P.[\zeta_{{\rm Bh}}(\alpha|\bm{w})]_{\alpha=q}=
\frac{(-1)^q}{(q-1)!}
\sum_{m=0}^{d-q}\frac{_dS_m^{(d)}(0)}{m!(d-q-m)!}
F\left[x^{d-q-m}\log x\right]_{x=\bm{w}}\\
\mbox{}
+\frac{_dS_{d-q}^{(d)}(0)(-1)^{d+q+1}}{(q-1)!(d-q)!}H_{q-1}
+\sum_{\bm{n}\in\mathbb{N}_0^d\setminus\{\bm{0}\}}\Bigg[
\frac{1}{(\bm{n}\cdot\bm{w})^{q}}
\\
\mbox{}\left.
+\frac{(-1)^q}{(q-1)!}\sum_{m=0}^{d-q}\frac{_dS_m^{(d)}(0)}{m!(d-q-m)!}G\left[(\bm{n}\cdot\bm{w}+x)^{d-q-m}\log(\bm{n}\cdot\bm{w}+x)\right]_{x=\bm{w}}\right]  \; .
%\label{fp_Bh_series}
\end{multline*}
This can be presented in the form of a limit as
\setlength\multlinegap{0pt}
\begin{multline*}
F.P.[\zeta_{{\rm Bh}}(\alpha|\bm{w})]_{\alpha=q}=
\lim_{M\tends\infty}%\Biggl[ 
\left[ 
\frac{(-1)^q}{(q-1)!}
%\\
%\mbox{}\times
\sum_{m=0}^{d-q}\frac{_dS_m^{(d)}(0)}{m!(d-q-m)!}
F\left[x^{d-q-m}\log x\right]_{x=\bm{w}}M^{d-q-m}
\right.
\\
\left. \mbox{}
+\frac{_dS_{d-q}^{(d)}(0)(-1)^{d+q+1}}{(q-1)!(d-q)!}\log M
%\\
%\mbox{}\left. 
+\sum_{\bm{n}\in\mathcal{C}_{M-1}\setminus\{\bm{0}\}}(\bm{n}\cdot\bm{w})^{-q}\right]
+\frac{_dS_{d-q}^{(d)}(0)(-1)^{d+q+1}}{(q-1)!(d-q)!}H_{q-1}\; .
%\label{fp_Bh_limit}
\end{multline*}
\setlength\multlinegap{10pt}
(We have used the fact that $F[x^l]_{x=\bm{w}}=0$ for $l=1,\ldots,d-1$ and $F[1]_{x=\bm{w}}=(-1)^{d-1}$.)

The derivative at $\alpha=0$ is given in series form as
\begin{multline*}
\zeta_{{\rm Bh}}'(0|\bm{w})=
\sum_{\bm{n}\in\mathbb{N}_0^d\setminus\{\bm{0}\}}\left[
-\log(\bm{n}\cdot\bm{w})-H_d
+\sum_{m=0}^{d}\frac{_dS_m^{(d)}(0)}{m!(d-m)!}
\right. \\
\mbox{}\times
G\left[ (\bm{n}\cdot\bm{w}+x)^{d-m}\log(\bm{n}\cdot\bm{w}+x)\right]_{x=\bm{w}}\Biggr] +
\sum_{m=0}^{d} \frac{_dS_m^{(d)}(0)}{m!(d-m)!}F\left[ x^{d-m}\log x\right]_{x=\bm{w}}-H_d
%\label{derivative_Bh_series}
\end{multline*}
and in limit form as
\begin{multline*}
\zeta_{{\rm Bh}}'(0|\bm{w})=\lim_{M\tends\infty}\left[ M^d(\log M-H_d)
+\sum_{m=0}^{d}\frac{_dS_m^{(d)}(0)}{m!(d-m)!}F\left[ x^{d-m}\log x\right]_{x=\bm{w}}M^{d-m}
\right. \\
\mbox{}\left.+(-1)^{d-1}\frac{_dS_d^{(d)}(0)}{d!}\log M -\sum_{\bm{n}\in\mathcal{C}_{M-1}\setminus\{\bm{0}\}}
\log(\bm{n}\cdot\bm{w})\right] \; .
\end{multline*}

Note that the finite parts of $\zeta_{{\rm Bh}}$ can be obtained directly from the respective finite parts of $\zeta_{{\rm B}}$ by removing the $1/a^q$ term and setting $a=0$. However, we present them here in a slightly more convenient form. Likewise, $\zeta_{{\rm Bh}}'(0|\bm{w})$ can be obtained from $\zeta_{{\rm B}}'(0,a|\bm{w})$ by simply removing the term $-\log a$ and setting $a=0$ in the remaining expression.

From these observations and Eqs.~(\ref{def_rho}), (\ref{def_gamma_dq}) and (\ref{relation_FP_Psi}), we have the additional relations
\begin{equation}
\zeta_{{\rm Bh}}'(0,\bm{w})=-\log \rho_{{\rm B}}(\bm{w})
\label{relation_rho_zetaBh}
\end{equation}
and
\begin{equation}
F.P.\left[ \zeta_{{\rm Bh}}(\alpha|\bm{w})\right]_{\alpha=q}=\frac{(-1)^{q-1}}{(q-1)!}\gamma_{dq}(\bm{w}) -H_{q-1}Res\left[\zeta_{{\rm Bh}}(\alpha|\bm{w})\right]_{\alpha=q}\; .
\label{relation_FP_gamma_dq}
\end{equation}
Thus, the series and limit representations we have obtained, used in conjunction with Eqs.~(\ref{def_Gamma}), (\ref{relation_FP_Psi}), (\ref{relation_rho_zetaBh}) and (\ref{relation_FP_gamma_dq}), immediately yield corresponding representations for $\log \rho_{{\rm B}}(\bm{w})$, $\log\Gamma_{{\rm B}}(a|\bm{w})$, $\Psi_{{\rm B}} ^{(q)}(a|\bm{w})$ and $\gamma_{dq}(\bm{w})$.

Before leaving this section, we observe that it is possible that other multiple $\zeta$-functions could have representations similar to the one in (\ref{Barnes_series1}). If so, this provides an alternative route to analytic continuation. It would be interesting to investigate this matter.

\subsection{Special case: $d=2$}

\label{serieslimits_2D}

The expressions we have obtained so far are substantially simplified when we consider specific cases. By way of illustration, we consider here the $d=2$ case. The simplification is large in both series form and limit form, although it is somewhat larger in the latter. For brevity, we present only the results in their limit form. We have then, in the $d=2$ case, from (\ref{fp_B_lim2}) and (\ref{derivative_B_lim}) and using in addition (\ref{Bernoulli}) and (\ref{Bernoullian}) for the computation of the needed Bernoullian functions,
\begin{multline*}
F.P.[\zeta_{{\rm B}}(\alpha,a|(w_1,w_2))]_{\alpha=2}=\lim_{M\tends\infty}\left[
-\frac{1}{w_1w_2}\log M+\sum_{\bm{n}\in\mathcal{C}_{M-1}}\frac{1}{(a+\bm{n}\cdot\bm{w})^2}\right]
\\
\mbox{}
+\frac{1}{w_1w_2}\left( -1+\log\frac{w_1+w_2}{w_1w_2}\right) \; ,\\
\end{multline*}
\begin{multline*}
F.P.[\zeta_{{\rm B}}(\alpha,a|(w_1,w_2))]_{\alpha=1}=
\lim_{M\tends\infty}\left[
-\left(\frac{1}{w_2}\log\frac{w_1+w_2}{w_1}+\leftrightarrow\right)M
\right.
\\ 
\left.\mbox{}
-\frac{1}{w_1w_2}\left(\frac{w_1+w_2}{2}-a\right)\log M 
%\right.
%\\ 
%\left.\mbox{}
+\sum_{\bm{n}\in\mathcal{C}_{M-1}}\frac{1}{a+\bm{n}\cdot\bm{w}}\right] \\
\mbox{}
+\frac{1}{w_1w_2}\left( \frac{w_1+w_2}{2}-a\right)\log\frac{w_1+w_2}{w_1w_2}
\; ,
\end{multline*}
\begin{multline*}
\zeta_{{\rm B}}'(0,a|\bm{w})=\lim_{M\tends\infty}\left[ M^2\log M+
\left(\frac{w_1}{2w_2}\log\frac{w_1+w_2}{w_1}+\leftrightarrow +\log(w_1+w_2)-\frac{3}{2}\right) M^2
\right. \\
\mbox{}+\left(\frac{2a-w_1-w_2}{2w_2}\log\frac{w_1+w_2}{w_1}+\leftrightarrow\right) M
\\
\mbox{}\left. -\frac{1}{2w_1w_2}\left( a^2-(w_1+w_2)a+\frac{(w_1+w_2)^2+w_1w_2}{6}\right)\log M
-\sum_{\bm{n}\in\mathcal{C}_{M-1}}\log(a+\bm{n}\cdot\bm{w})\right] 
\\
\mbox{}
+\frac{1}{2w_1w_2}\left[ a^2-(w_1+w_2)a+\frac{(w_1+w_2)^2+w_1w_2}{6}\right] \log\frac{w_1+w_2}{w_1w_2}\, .
\end{multline*}
The meaning of the symbol $\leftrightarrow$ is that of symmetrization. It stands for a term which is identical to the one immediately preceding it with the roles of $w_1$ and $w_2$ exchanged. Naturally, in the present case $\mathcal{C}_{M-1}=\{0,\ldots,M-1\}^2$, since $d=2$.

From these results, the ones for $\zeta_{{\rm Bh}}$ follow in a straightforward fashion: we merely omit the $\bm{n}=(0,0)$ term from the $\bm{n}\in\{ 0,1,\ldots,M-1\}^2$ sum and set $a=0$.

The above expressions, as well as those pertaining to the $d=3$ case (which, for reasons of space we do not display), have been confirmed numerically. For this, we checked them against results coming from a numerical treatment of the integral representations of the next section and, in a few cases where it is possible and simple ($\bm{w}=(1,n)$ in $d=2$ and $\bm{w}=(1,n,n)$, $\bm{w}=(1,1,n)$ in $d=3$, $n\in\mathbb{N}$), against the expressions which result from casting $\zeta_{{\rm B}}$ as a finite combination of Hurwitz $\zeta$-functions.

\section{Integral representations}

\label{integrals}

In his work, Barnes gave contour and line integral representations for $\log\Gamma_{{\rm B}}$, $\log\rho_{{\rm B}}$, $\Psi_{{\rm B}} ^{(q)}$ and $\gamma_{dq}$\footnote{\label{rho_error}Barnes' line integral representation for $\log\rho_{{\rm B}}(\bm{w})$, given in p.~411 of \cite{Barnes1904a}, contains an error: there should be an extra term, $(-1)^{r-1} \mbox{}_rS_1'(0){\rm e}^{-z{\rm e}^{-{\rm i}\phi}}$, inside the curly brackets in the integrand, and an extra term $(-1)^{r-1} \mbox{}_rS_1'(0){\rm i}\phi$ outside the integral. This error can be traced to the very last step in the derivation of the integral and is due to the oversight of a term $_rS_1'(0)$.}. Through Eqs.~(\ref{def_Gamma}), (\ref{relation_FP_Psi}), (\ref{relation_rho_zetaBh}) and (\ref{relation_FP_gamma_dq}), these representations immediately provide analogous representations for the finite parts and derivative at zero of $\zeta_{{\rm B}}$ and $\zeta_{{\rm Bh}}$. In this section, we provide alternative line integral representations 
for the finite parts of $\zeta_{{\rm B}}$ and $\zeta_{{\rm Bh}}$ and the derivative at zero of $\zeta_{{\rm Bh}}$ (a related expression for the derivative at zero of $\zeta_{{\rm B}}$ was given by Ruijsenaars in Eq.~(3.13) of \cite{Ruijsenaars2000}).
The representations here obtained, given below in Eqs.~(\ref{fp_B_integral}), (\ref{fp_Bh_integral2}) and (\ref{derivative_Bh_integral2}) (together with Eq.~(3.13) of \cite{Ruijsenaars2000}), in turn provide line integral representations for the Barnes' functions just mentioned, alternative to the ones given by Barnes.

In contrast to the series and limit representations of the previous section, whose validity is absolutely general, the integral representations of the present section are valid when $\Re(a)>0$ and $\Re(w_i)>0$ (corresponding to $\theta=0$ in the definition of the half-plane $\mathcal{H}$).

\subsection{Barnes ${\bm\zeta}$-function}

\label{integrals_B}

Integral representations for the finite parts of $\zeta_{{\rm B}}$ are available in a very direct way from an integral representation for $\zeta_{{\rm B}}$ given by Ruijsenaars (Eq.~(3.8) in \cite{Ruijsenaars2000}). This representation can be cast as
\begin{multline}
\zeta_{{\rm B}}(\alpha,a|\bm{w})=\frac{1}{\Gamma (\alpha)\prod_{i=1}^dw_i}\sum_{k=0}^M\frac{(-1)^kB_k(\bm{w})\Gamma(\alpha -d+k)}{k!a^{\alpha -d+k}}+I_M(\alpha,a|\bm{w})\, ,\\ \Re(\alpha)>d-M-1 
\; ,
\label{B_Ruijsenaars}
\end{multline}
where
\begin{equation*}
I_M(\alpha,a|\bm{w})=\frac{1}{\Gamma(\alpha)}\int_0^{\infty}{\rm d}t\, {\rm e}^{-at}\left(
\frac{t^{\alpha -1}}{\prod_{i=1}^d(1-{\rm e}^{-w_it})}-\frac{t^{\alpha-d-1}}{\prod_{i=1}^dw_i}
\sum_{k=0}^M(-1)^kB_k(\bm{w})\frac{t^k}{k!}\right)
\end{equation*}
and $B_k(\bm{w})$ are the higher order Bernoulli numbers, defined from the higher order Bernoulli polynomials by $B_k(\bm{w})=B_k(0|\bm{w})$. (As mentioned in footnote \ref{Bernoulli_def}, Ruijsenaars uses a slightly different convention for these objects.) $M$ can be chosen at will from $0,1,2,\ldots$ 
A specific case of this representation had been obtained before by Holthaus et al. \cite{HolthausKalinowskiKirsten1998}.

The setting in \cite{Ruijsenaars2000} assumes $\Re (a)>0$ and $w_1,\ldots,w_d\in\mathbb{R}^+$. However, Eq.~(\ref{B_Ruijsenaars}) is valid under the more general conditions $\Re (a)>0$, $\Re (w_i)>0$, $i=1,\ldots,d$. Indeed, the line integral representation for $\zeta_{{\rm B}}$ which is used as a starting point in the derivation of Eq.~(\ref{B_Ruijsenaars}) is valid in this enlarged parameter range. This can be seen, for example, from the contour integral representation for $\zeta_{{\rm B}}$, by the standard procedure of reducing the Hankel contour to a line.

The integral $I_M(\alpha,a|\bm{w})$ is an analytic function of $\alpha$ for $\Re(\alpha)>d-M-1$. Therefore, the pole structure of $\zeta_{{\rm B}}$ is contained in the $\Gamma$-functions of the first term in (\ref{B_Ruijsenaars}). To obtain the finite part at $\alpha=q$ we merely compute the finite part of the first term, from knowledge of the Laurent expansions of the $\Gamma$-functions, and set $\alpha=q$ in $I_M(\alpha,a|\bm{w})$, taking $M=d-q$. Thus,
\begin{multline}
F.P.\left[\zeta_{{\rm B}}(\alpha,a|\bm{w})\right]_{\alpha=q}=
I_{d-q}(q,a|\bm{w})
 \\
\mbox{}+\frac{(-1)^{d-q}}{\left(\prod_{i=1}^dw_i\right)(q-1)!}
\sum_{k=0}^{d-q}\frac{a^{d-q-k}B_k(\bm{w})}{k!(d-q-k)!}(H_{d-q-k}-H_{q-1}-\log a)\; .
\label{fp_B_integral}
\end{multline}

A similar integral representation for the derivative of $\zeta_{{\rm B}}$ at $\alpha=0$ is obtained from (\ref{B_Ruijsenaars}) by differentiation (Eq.~(3.13) of \cite{Ruijsenaars2000}).

\subsection{Homogeneous Barnes ${\bm\zeta}$-function}

\label{integrals_Bh}

To treat the homogeneous case, we could think of subtracting a term $1/a^{\alpha}$ from the representation (\ref{B_Ruijsenaars}) and taking $a\tends 0$. However, both terms of (\ref{B_Ruijsenaars}) become divergent in this limit. Another possibility would be to cast $\zeta_{{\rm Bh}}$ as a combination of several $\zeta_{{\rm B}}$ of dimensions $1,2,\ldots,d$ and use (\ref{B_Ruijsenaars}) to obtain the finite parts as in section~\ref{integrals_B}. This was done in \cite{HolthausKalinowskiKirsten1998} for a finite part in the $d=2$ setting, in which case the procedure is quite convenient.
It can be applied to the general $d$-dimensional case, although the resulting expressions become substantially more complicated, particularly when dealing with the poles towards the left or the derivative at $\alpha=0$. We give here a streamlined procedure which treats $\zeta_{{\rm Bh}}$ directly, 
whereby the analytic continuation of $\zeta_{{\rm Bh}}$ is obtained in a form akin to (\ref{B_Ruijsenaars}),
thus yielding a considerably simpler end result in the general case. It is based on the same idea of adding and subtracting of divergencies as used in \cite{HolthausKalinowskiKirsten1998,Ruijsenaars2000}, with one additional feature that enables its application to the homogeneous case.

Using the $\Gamma$-function identity
\begin{equation*}
(\bm{n}\cdot\bm{w})^{-\alpha}\,\Gamma(\alpha)=\int_0^\infty t^{\alpha -1}{\rm e}^{-(\bm{n}\cdot\bm{w})t}\, {\rm d}t\, ,
\quad \Re(\alpha)>0\, ,\; \Re(\bm{n}\cdot\bm{w})>0\, ,
\end{equation*}
in (\ref{zetaBh}), we have
\begin{equation*}
\begin{split}
\Gamma(\alpha)\zeta_{{\rm Bh}}(\alpha|\bm{w})&=
\lim_{M\tends\infty}\sum_{\bm{n}\in\mathcal{C}_M\setminus\{\bm{0}\}}\int_0^\infty t^{\alpha -1}{\rm e}^{-(\bm{n}\cdot\bm{w})t}\, {\rm d}t
\\&
=\lim_{M\tends\infty}\int_0^\infty t^{\alpha -1}\left(
\prod_{i=1}^d\frac{1-{\rm e}^{-(M+1)w_it}}{1-{\rm e}^{-w_it}}-1\right) \, {\rm d}t
\\&\hspace{-65pt}
=\lim_{M\tends\infty}\int_0^\infty t^{\alpha -1}\left(
\prod_{i=1}^d\frac{1}{1-{\rm e}^{-w_it}}-1+\text{terms of the form }\frac{{\rm e}^{-(M+1)wt}}{\prod_{i=1}^d(1-{\rm e}^{-w_it})}
\right) \, {\rm d}t\; ,
\end{split}
\end{equation*}
where $w$ is a finite sum of the $w_i$'s. (We used limits in order to justify the interchange of sum and integral in the second line above.) Now, when $M\tends\infty$ the last part of the integral vanishes and we are left with
\begin{equation}
\Gamma(\alpha)\zeta_{{\rm Bh}}(\alpha,\bm{w})=\int_0^\infty t^{\alpha -1}\left(
\prod_{i=1}^d\frac{1}{1-{\rm e}^{-w_it}}-1\right) \, {\rm d}t\; ,\quad
\Re(\alpha)>d\, ,\; \Re(w_i)>0\, .
\label{line_zetaBh}
\end{equation}

Let $c$ be any constant satisfying $\Re(c)>0$ and let $M\in\mathbb{N}_0$. We multiply and divide the integrand in (\ref{line_zetaBh}) by ${\rm e}^{ct}$. This enables us to add and subtract the divergencies contained in the lower limit of integration when $\Re(\alpha)\leq d$, using 
the first terms of the expansion of ${\rm e}^{ct}/\prod_{i=1}^d(1-{\rm e}^{-w_it})$ from (\ref{Bernoulli}),
and split the integral, in the following way.
\begin{multline}
\zeta_{{\rm Bh}}(\alpha|\bm{w})
=\frac{1}{\Gamma(\alpha)}\int_0^{\infty}{\rm d}t\, t^{\alpha -1}{\rm e}^{-ct}
\left( \frac{{\rm e}^{ct}}{\prod_{i=1}^d(1-{\rm e}^{-w_it})}-{\rm e}^{ct}\right)\\
=\frac{1}{\Gamma(\alpha)\prod_{i=1}^dw_i}\sum_{k=0}^M\frac{(-1)^kB_k(-c|\bm{w})\Gamma(\alpha -d+k)}{k!c^{\alpha -d+k}}\\
-\frac{1}{\Gamma(\alpha)c^{\alpha}}\sum_{k=0}^{M-d}\frac{\Gamma(\alpha+k)}{k!}
+J_{M,c}(\alpha|\bm{w}) \; ,
\label{zetaBh_integral}
\end{multline}
where
\begin{multline}
J_{M,c}(\alpha|\bm{w})=\frac{1}{\Gamma(\alpha)}\int_0^\infty {\rm d}t\, t^{\alpha -1}\left(
\frac{1}{\prod_{i=1}^d(1-{\rm e}^{-w_it})}-1\right. \\
\left. 
\mbox{}-\frac{t^{-d}{\rm e}^{-ct}}{\prod_{i=1}^dw_i}\sum_{k=0}^M
(-1)^kB_k(-c|\bm{w})\frac{t^k}{k!}
+{\rm e}^{-ct}\sum_{k=0}^{M-d}\frac{(ct)^k}{k!}
\right) \; .
\label{Jintegral_zetaBh}
\end{multline}
By convention, empty sums (the $k=0$ to $M-d$ sums in (\ref{zetaBh_integral}) and (\ref{Jintegral_zetaBh}) if $M<d$) have value zero.
This representation can be regarded as the analogue of (\ref{B_Ruijsenaars}) for the homogeneous case. It provides the analytic continuation of $\zeta_{{\rm Bh}}$ to $\Re(\alpha)>d-M-1$ and it is valid for $\Re(w_i)>0$, $i=1,\ldots,d$. A convenient choice for the constant $c$ is $c=1$. 

We obtain the finite parts as before. Setting $c=1$, we have
\begin{multline}
F.P.\left[\zeta_{{\rm Bh}}(\alpha|\bm{w})\right]_{\alpha=q}=J_{d-q,1}(q|\bm{w})\\
\mbox{}+\frac{(-1)^{d-q}}{\left(\prod_{i=1}^dw_i\right)(q-1)!}
\sum_{k=0}^{d-q}\frac{B_k(-1|\bm{w})}{k!(d-q-k)!}(H_{d-q-k}-H_{q-1}) \; .
\label{fp_Bh_integral}
\end{multline}

In order to obtain $\zeta_{{\rm Bh}}'(0|\bm{w})$, we set $M=d$ and differentiate, which yields
\begin{multline}
\zeta_{{\rm Bh}}'(0|\bm{w})=
\frac{(-1)^{d}}{\prod_{i=1}^dw_i}\sum_{k=0}^{d-1}\frac{B_k(-1|\bm{w})H_{d-k}}{k!(d-k)!}\\
\mbox{}+\int_0^\infty {\rm d}t\, t^{-1}\left(
\frac{1}{\prod_{i=1}^d(1-{\rm e}^{-w_it})}-1+{\rm e}^{-t}-\frac{t^{-d}{\rm e}^{-t}}{\prod_{i=1}^dw_i}
\sum_{k=0}^d(-1)^kB_k(-1|\bm{w})\frac{t^k}{k!}\right)\, .
\label{derivative_Bh_integral}
\end{multline}

These can be simplified by using $B_k(x|\bm{w})=\sum_{l=0}^k\binom{k}{l}x^lB_{k-l}(\bm{w})$ \cite{Norlund1922}. For this purpose, note that
\begin{equation}
\begin{split}
\sum_{k=0}^{d-q}\frac{B_k(-1|\bm{w})}{k!(d-q-k)!}&=\sum_{k=0}^{d-q}\sum_{l=0}^k\frac{(-1)^lB_{k-l}(\bm{w})}{(d-q-k)!l!(k-l)!}\\
&=\sum_{j=0}^{d-q}\frac{B_j(\bm{w})}{j!}\sum_{l=0}^{d-q-j}\frac{(-1)^l}{(d-q-j-l)!l!}=\frac{B_{d-q}(\bm{w})}{(d-q)!}
\end{split}
\label{binom1}
\end{equation}
and
\begin{equation}
\begin{split}
\sum_{k=0}^{d-q}\frac{B_k(-1|\bm{w})}{k!(d-q-k)!}H_{d-q-k}&=\sum_{k=0}^{d-q-1}\frac{H_{d-q-k}}{k!(d-q-k)!}\sum_{l=0}^k\binom{k}{l}(-1)^lB_{k-l}(\bm{w})\\
&=\sum_{j=0}^{d-q-1}\frac{B_j(\bm{w})}{j!}\sum_{l=0}^{d-q-1-j}
\frac{(-1)^lH_{d-q-j-l}}{(d-q-j-l)!l!} \; .
\end{split}
\label{binom2}
\end{equation}
The inner sum in the last expression can be simplified further. Let $m=d-q-j$. We have
\begin{equation}
\begin{split}
\sum_{l=0}^{m-1}\frac{(-1)^lH_{m-l}}{(m-l)!l!}&=\sum_{k=1}^m\sum_{l=0}^{m-k}\frac{(-1)^l}{(m-l)!l!k}
=\frac{1}{m!}\sum_{k=1}^m\frac{(-1)^{m-k}}{k}\binom{m-1}{m-k}\\ &
=\frac{1}{m!m}\sum_{k=1}^m(-1)^{m-k}\binom{m}{k}=\frac{(-1)^{m+1}}{m!m} \; .
\end{split}
\label{binom3}
\end{equation}
Using (\ref{binom1}), (\ref{binom2}) and (\ref{binom3}) in (\ref{fp_Bh_integral}) and (\ref{derivative_Bh_integral}), we obtain
\begin{multline}
F.P.\left[\zeta_{{\rm Bh}}(\alpha|\bm{w})\right]_{\alpha=q}=\frac{1}{\left(\prod_{i=1}^dw_i\right)(q-1)!}
\Biggl[ 
\frac{(-1)^{d-q+1}B_{d-q}(\bm{w})H_{q-1}}{(d-q)!}
\\ \left. 
+\sum_{j=0}^{d-q-1}\frac{(-1)^{j+1}B_j(\bm{w})}{j!(d-q-j)!(d-q-j)}\right] 
+J_{d-q,1}(q|\bm{w})
\label{fp_Bh_integral2}
\end{multline}
and
\begin{multline}
\zeta_{{\rm Bh}}'(0|\bm{w})=\frac{1}{\prod_{i=1}^dw_i}\sum_{j=0}^{d-1}\frac{(-1)^{j+1}B_j(\bm{w})}{j!(d-j)!(d-j)}\\
+\int_0^\infty {\rm d}t\, t^{-1}\left(
\frac{1}{\prod_{i=1}^d(1-{\rm e}^{-w_it})}-1+{\rm e}^{-t}-\frac{t^{-d}{\rm e}^{-t}}{\prod_{i=1}^dw_i}
\sum_{k=0}^d(-1)^kB_k(-1|\bm{w})\frac{t^k}{k!}\right)\, .
\label{derivative_Bh_integral2}
\end{multline}

The representations (\ref{fp_B_integral}), (\ref{fp_Bh_integral2}) and (\ref{derivative_Bh_integral2}) bear resemblances to Barnes' line integral representations for $\log\Gamma_{{\rm B}}$, $\log\rho_{{\rm B}}$, $\Psi_{{\rm B}} ^{(q)}$ and $\gamma_{dq}$ mentioned above. In fact, it is possible to show directly, starting from the integrals, the equivalence between both forms of representation (see however footnote~\ref{rho_error}).

\end{document}